\documentclass[a4paper,10pt]{article}
\usepackage{amsmath,amssymb,amsthm}

\theoremstyle{plain}
\newtheorem{thm}{Theorem}
\newtheorem{lem}[thm]{Lemma}
\newtheorem{prop}[thm]{Proposition}
\newtheorem{conj}[thm]{Conjecture}

\theoremstyle{definition}
\newtheorem{defn}{Definition}

\title{Constraints on counterexamples to the Casas-Alvero conjecture, and a verification in degree $12$}
\author{Wouter Castryck, Robert Laterveer, Myriam Ouna\"ies}

\begin{document}

\maketitle
\date{}

\begin{abstract}

  \noindent In a first (theoretical) part of this paper, we prove a number of constraints on
  hypothetical counterexamples to the Casas-Alvero conjecture,
  building on ideas of Graf von Bothmer, Labs, Schicho and \mbox{van de Woestijne}
  that were recently reinterpreted by Draisma and de Jong in terms of $p$-adic
  valuations. In a second (computational) part, we present ideas improving upon
  Diaz-Toca and Gonzalez-Vega's
  Gr\"obner basis approach to the Casas-Alvero conjecture.
  One application is
  an extension of the proof of Graf von Bothmer et al.\
  to the cases $5p^k$, $6p^k$ and $7p^k$ (that is, for each of these cases,
  we elaborate the finite list of primes $p$ for which their proof is not applicable).
  Finally, by combining both parts, we
  settle the Casas-Alvero conjecture in degree $12$ (the smallest open case).
\end{abstract}

\tableofcontents
\vspace{0.3cm}
\noindent \textbf{Files:} \verb"CAbadprimes.m" \ \verb"CAbadprimes7test.m" \ \verb"badprimes7.txt" \ \verb"CAdeg12.m"

\section{Introduction and overview} \label{introduction}

\textbf{(\ref{introduction}.1)} The subject of this article is the following intriguing conjecture \cite{CA}:
\begin{conj}[The Casas-Alvero conjecture, 2001] \label{CAconjecture} Let $f(x) \in \mathbb{C}[x]$ be of degree $d >0$
and suppose that for each $j=1, \dots, d-1$ there exists an $a \in \mathbb{C}$
such that $f(a) = f^{(j)}(a) = 0$, where $f^{(j)}(x)$ denotes the $j$th derivative.
Then $f(x)$ is the $d$th power of a linear polynomial.
\end{conj}
\noindent For each given degree $d$, proving Conjecture~\ref{CAconjecture} (if true) boils
down to a finite Gr\"obner basis computation. In 2006, this was used
by Diaz-Toca and Gonzalez-Vega to verify the conjecture
for $d \leq 7$ \cite{DTGV}. Shortly after, Graf von Bothmer, Labs,
Schicho and van de Woestijne \cite{GLSW} proved a theoretical result settling the
cases $d=p^k$ and $d=2p^k$ (where $p$ is prime and $k \geq 0$ is an integer).
The proof uses reduction-mod-$p$ arguments in algebraic geometry. It was
recently rewritten in the more elementary (and slightly more powerful) language of $p$-adic valuations, in a nice overview
due to Draisma and de Jong \cite{DD}.\\

\noindent \textbf{(\ref{introduction}.2)} By lack of a general strategy, beyond the degree, we subdivide the set of hypothetical counterexamples $f(x)$ to
the Casas-Alvero conjecture by
\begin{itemize}
 \item their number of distinct roots $\#\text{roots}(f)$,
 \item their \emph{type} $\text{type}(f)$, which is the minimal number of recycled roots minus one
 \[ \min \left\{ \, \#S \, \left| \, S \subset \mathbb{C} \text{ and } \forall j :
\exists \, a \in S : f(a) = f^{(j)}(a) = 0 \right. \, \right\} \, - \, 1 \]
 where $j$ ranges over $\{1, \dots, d-1\}$,
 \item their \emph{scenario} $\text{scen}(f)$, which is
\begin{equation} \label{defscen}
 \min \left\{ \, \left. (s_1, \dots, s_{d-1}) \in \mathbb{Z}_{\geq 0}^{d-1} \, \right| \,
 \exists \, \text{$a_i$'s} \in \mathbb{C}: \forall j : f(a_{s_j}) = f^{(j)}(a_{s_j}) = 0 \, \right\}
\end{equation}
 where the minimum is taken lexicographically and
 $j$ ranges over $\{1, \dots, d-1\}$. Note that $\text{type}(f)$ is the maximal entry of $\text{scen}(f)$.\\
\end{itemize}

\noindent \textbf{(\ref{introduction}.3)} The scenario $(s_1, \dots, s_{d-1}) \in \mathbb{Z}_{\geq 0}^{d-1}$ of
a degree $d$ counterexample $f \in \mathbb{C}[x]$ to the Casas-Alvero conjecture always satisfies
$s_1 = 0$ and $s_j \leq \max \{ \, s_i \, | \, i < j \, \}  + 1$
for all $j = 2, \dots, d-1$. A sequence of this form will therefore be called \emph{a scenario for degree $d$}.
In view of the above, the \emph{type} of a scenario is defined to be its maximal entry --
we denote it by $\text{type}(s)$.
The number of scenarios for a given degree $d$ grows quickly with $d$. E.g., in our main case
of interest $d=12$, we have
\[ 1, 1023, 28501, 145750, 246730, 179487, 63987, 11880, 1155, 55, 1 \]
scenarios of type $0, \dots, 10$, respectively, amounting to a total of $678570$.\\

\noindent \textbf{(\ref{introduction}.4)} Let $s = (s_1, \dots, s_{d-1})$ be a scenario for degree $d$, and let $t = \text{type}(s)$.
Let $f(x) \in \mathbb{C}[x]$ be a
degree $d$ counterexample to the Casas-Alvero conjecture. Then we say that $f(x)$ \emph{matches with $s$}
if there exist $a_0, \dots, a_t \in \mathbb{C}$ such that
\begin{itemize}
\item $f(x) = g(x) \cdot (x-a_0)(x-a_1)\cdots(x-a_t)$ for a degree $d-1-t$ polynomial $g(x) \in \mathbb{C}[x]$,
\item $f(a_{s_j}) = f^{(j)}(a_{s_j}) = 0$ for all $j=1, \dots, d-1$.
\end{itemize}
Clearly $f(x)$ matches with its own scenario $\text{scen}(f)$, but it may also match with various other scenarios.\\

\noindent \emph{Example.} Since it is conjecturally impossible to give examples over $\mathbb{C}$, consider
$f(x) = x(x-1)^4(x-8)(x-18) \in \mathbb{F}_{23}[x]$. One checks that the common roots of $f$ with $f^{(1)}, \dots, f^{(6)}$
are
\[ \{1\}, \quad \{1,18\}, \quad \{1\}, \quad \{0\}, \quad \{18\}, \quad \{1\},\]
respectively. So $\text{type}(f) = 2$ and $\text{scen}(f) = (0,0,0,1,2,0)$ (take $a_0=1,a_1=0,a_2=18$).
However, $f(x)$ also matches with $(0,1,0,2,1,0)$ (and many more).\\

\noindent \textbf{(\ref{introduction}.5)} In Section~\ref{general}, we prove a number of general constraints on these attributes. E.g.,
we find that
\begin{itemize}
 \item $\#\text{roots}(f) \geq 5$,
 \item $2 \leq \text{type}(f) \leq d-3$ (the first inequality being due to Draisma and Knopper \cite[Proposition~6]{DD}),
 \item if $\text{type}(f) = d-3$, then no consecutive entries of $\text{scen}(f)$ are equal.
\end{itemize}
The methods used here are classically flavoured (Gauss--Lucas, Newton, Rolle).\\

\noindent \textbf{(\ref{introduction}.6)} In Section~\ref{additional}, using the $p$-adic valuation approach,
we prove additional constraints for
certain special degrees.
Our main results are on degrees of the form $p+1$:

\begin{thm} \label{DeltaThm}
 Let $p$ be prime and let $f(x)$ be a degree $d=p+1$ counterexample to the Casas-Alvero conjecture.
 Let $c$ be the root of $f^{(d-1)}(x)$. Then $f^{(1)}(c) \neq 0$, and there exist at least two
 indices $2 \leq j_1 < j_2 \leq d-2$ such that $f^{(j_1)}(c) = f^{(j_2)}(c) = 0$.
 In particular, $\emph{type}(f) \leq d-4$.
 Moreover, if $j_1 < \dots < j_m$ are the indices between $2$ and $d-2$ for which
 $f^{(d - j_1)}(c) = \dots = f^{(d - j_m)}(c) = 0$, then the determinant of
 \begin{equation}\label{Delta}
\Delta_f=\left[
\begin{array}{c c c  c  c c }
-1& j_1 & 0  & 0 &\cdots  & 0\\
-1 & \binom{j_2-2}{j_1-2}j_2 & j_2  & 0 &\cdots & 0 \\
\vdots & \vdots & \vdots  &\vdots & & \vdots \\
-1 & \binom{j_m-2}{j_1-2}j_m & \binom{j_m-2}{j_{2}-2}j_m & &\cdots & j_m \\
-1 & (-1)^{j_1} & (-1)^{j_2} & & \cdots  & (-1)^{j_m}
\end{array}
\right]
\end{equation}
 is a multiple of $p$.
\end{thm}

\noindent Theorem~\ref{DeltaThm}
implies that every degree $d=p+1$ counterexample
to the Casas-Alvero conjecture matches with an element
of the strongly reduced list of scenarios
$s = (s_1,\dots, s_{d-1})$ for which
\begin{itemize}
  \item $s_{d-1} \ne 0$,
  \item the set of indices $2 \leq j \leq d-2$
for which $s_{d-j} = s_{d-1}$ satisfies the above determinant condition.
\end{itemize}
For $d=12$ ($p=11$), the
list contains
\begin{equation} \label{restrictedlist}
 0, 48, 1668, 8172, 11586, 6298, 1469, 146, 5, 0, 0
\end{equation}
scenarios of type $0, \dots, 10$, respectively, amounting
to a total of $29392$. In type $8$,
the five scenarios read
\begin{equation} \label{scensoftype8left}
 \begin{array}{l} (0,1,2,3,4,5,6,7,3,8,3), \\
                  (0,1,2,3,4,5,5,6,7,8,5), \\
                  (0,1,2,3,4,3,5,6,7,8,3), \\
                  (0,1,2,3,4,2,5,6,7,8,2), \\
                  (0,1,2,3,2,4,5,6,7,8,2); \\
 \end{array}
\end{equation}
indeed, the only pairs $(j_1,j_2)$ for which $\Delta_f \equiv 0 \bmod 11$ are $(3,8)$, $(5,6)$,
$(6,8)$, $(6,9)$, $(7,9)$.\\

\noindent \textbf{(\ref{introduction}.7)} For the computational part of
our paper, we turn back to the original reduction-mod-$p$ setting used
by Graf von Bothmer et al. Because of the interplay between characteristic $0$ and characteristic $p>0$, the following
general definition is convenient.
\begin{defn}
Let $k$ be an algebraically closed field. We say that a degree $d$ polynomial
$f \in k[x]$ ($d >0$) is a \emph{Casas-Alvero polynomial} or \emph{CA-polynomial} (over $k$)
if $f$ is not a power of a linear polynomial and if for each $j = 1, \dots, d - 1$ there exists an $a \in k$ such that
$f(a) = f^{(j)}_H(a) = 0$.
\end{defn}
\noindent Here, $f^{(j)}_H$ denotes the $j$th Hasse derivative (using Hasse derivatives
turns the Casas-Alvero condition somewhat more restrictive -- it makes no difference in characteristic $0$ or
$p > d-1$, where $f^{(j)}_H = \frac{1}{j!} f^{(j)}$). Then the main theorem of \cite{GLSW} reads:

\begin{thm}[Graf von Bothmer, Labs, Schicho, van de Woestijne] \label{Grafetal}
 Let $d > 0$ be an integer and let $p$ be a prime number. If no CA-polynomials
of degree $d$ exist over $\overline{\mathbb{F}}_p$, then the Casas-Alvero conjecture
is true in degree $dp^k$ for all integers $k \geq 0$.
\end{thm}
\noindent Since it is trivial that no CA-polynomials of degree $1$ or $2$ can exist (in any characteristic),
the cases $p^k$ and $2p^k$ follow. More generally, we call a prime $p$ a \emph{bad prime for degree $d$}
if there exist CA-polynomials of degree $d$ in characteristic $p$. Then it is easily verified
that $p=2$ is the sole bad
prime for degree $d=3$. De Jong and Draisma \cite{DD} proved that the bad primes for degree $d=4$ are $p=3,5,7$.\\

\noindent \textbf{(\ref{introduction}.8)} In Section~\ref{computational}
we present an algorithm, the basic version of which takes
as input an integer $d > 0$ and
a prime number $p$ (or $p = 0$), and outputs whether or not CA-polynomials
of degree $d$ exist in characteristic $p$.
The basic idea is to classify all CA-polynomials by their scenario
(the definitions in \textbf{(\ref{introduction}.2)} straightforwardly generalize
to arbitrary $k$ -- this was already used in the example in \textbf{(\ref{introduction}.4)} there under). 
We will see that scenarios of
moderately low type $t$ can be ruled out easily (if the Casas-Alvero conjecture is true).
In characteristic $0$, the computation is feasible up to $d \cdot t \approx 50$, say. In moderate characteristic $p$,
this can be pushed to about twice that value.\\

\noindent \textbf{(\ref{introduction}.9)} By running the algorithm in characteristic $0$ and analyzing
the prime factors appearing in certain resulting Nullstellensatz expansions, we can find the bad
primes for $d$ up to $7$.
\begin{thm} \label{badprimesoverview} There are
\begin{itemize}
  \item $9$ bad primes for degree $d=5$, namely,
\[ p = 2, 3, 7, 11, 131, 193, 599, 3541, \text{and } 8009,\]
  \item $53$ bad primes for degree $d=6$, namely, the primes listed in Table~\ref{badprimes6},
  \item $366$ bad primes for degree $d=7$, namely, the primes listed in the file \verb"badprimes7.txt"
  that accompanies this paper -- the smallest non-bad prime (apart from $p=7$) is $127$ -- the largest bad prime
  is
\begin{center}
\footnotesize{$ \begin{array}{l} 24984712021698392647916525667237483011737174983678606896870094983849 \\ 9096141806825287856933123954724798488422551659890912229726792102063 \\ \end{array}$}
\end{center}
  (a $135$-digit number).
\end{itemize}
\end{thm}
\begin{table}
\begin{center}
\begin{tabular}{|l|l|l|l|}
\hline
  2 & 5 & 7 & 11 \\
13 & 19 & 23 & 29 \\
37 & 47 & 61 & 67 \\
73 & 97 & 257 & 811 \\
983 & 1069 & 1087 & 1187 \\
1487 & 1499 & 1901 & 2287 \\
3209 & 3877 & 3881 & 4019 \\
4943 & 5471 & 6983 & 8699 \\
9337 & 15131 & 15823 & 20771 \\
21379 & 23993 & 150203 & 266587 \\
547061 & 685177 & 885061 & 1030951 \\
7783207 & 17250187 & 40362599 & 9348983563 \\
70016757407 & 2610767527031 & 225833117528659 & 7390044713023799 \\
51313000813080529 &  &  & \\
\hline
\end{tabular}
\caption{Bad primes for degree $6$ ($53$ primes)}
\label{badprimes6}
\end{center}
\end{table}
\noindent We note that the bad primes for $d=5$ have been independently elaborated (by hand) by Chellali and Salinier \cite{CS}.\\

\noindent \textbf{(\ref{introduction}.10)} Finally in Section~\ref{degree12}, we combine our
theoretical and computational approaches. Naively running our algorithm in
degree $12$ lies completely out of reach.
But in view of Theorem~\ref{DeltaThm} and certain reduction-mod-$p$ considerations,
it suffices to restrict the algorithm to a limited list of scenarios, and to run it in characteristic $p$.
As such, the computation becomes feasible:
\begin{thm} \label{deg12thm}
 Conjecture~\ref{CAconjecture} is true for $d=12$.
\end{thm}
\noindent The margin is tight: each of the five scenarios of (\ref{scensoftype8left}) took approximately
three weeks of computation and required about $90$ GB of
RAM. Pushing the analogous computation to $d=20$, the next open case, is utopic.\\

\noindent \textbf{(\ref{introduction}.11)}  The main computations have been carried out using Magma \cite{magma} version 2.18-2 on
a computer called \verb"matrix", running
Ubuntu 11.10 on a $6$-core Intel Xeon 2.53 GHz processor with 96 GB RAM. Some additional calculations were executed
using Magma version 2.15-12 on \verb"kasparov", running Debian GNU/Linux 6.0.4 on an $8$-core x86-64 2.93 GHz processor with 64 GB RAM.\\

\noindent \textbf{(\ref{introduction}.12)} We would like to thank Filip Cools, Jan Schepers and Fr\'ederik Vercauteren
for some helpful discussions. We are also grateful to the Department of Electrical Engineering (KU Leuven), for
allowing us to use \verb"kasparov".

\section{General constraints on counterexamples} \label{general}

\noindent \textbf{(\ref{general}.1)}
The following easy fact will be used throughout:
\begin{lem}\label{change}
Let $f$ be a CA-polynomial over $k$ of degree $d > 0$,
$\alpha_1, \alpha_2 \in k^*$ and $\beta\in k$. Then
the polynomial $g(x)=\alpha_1 f(\alpha_2 x+\beta)$ is also CA.
\end{lem}
\noindent The polynomials $f$ and $g$ will be called \emph{equivalent}. Note that
the number of distinct roots, the type, the scenario, the matching or not with a given scenario, \dots \ are all preserved by equivalence.\\

\noindent \textbf{(\ref{general}.2)} We begin with some
considerations on the type:

\begin{prop} \label{gonality}
Let $f \in \mathbb{C}[x]$ be a CA-polynomial of degree $d$ and let
$\Gamma$ be the convex hull of the roots of $f$ (when plotted in the complex plane). Let
$m \geq 2$ be the maximum of the multiplicities of these roots, and
let $\delta = 1$ if this maximum is attained by a non-vertex of
$\Gamma$ (let $\delta = 0$ otherwise). Let
$\gamma \geq 2$ be the number of vertices
of $\Gamma$.
Then $2 \leq \emph{type}(f) \leq d + 1 - \gamma - m - \delta \leq d-3$.
\end{prop}

\noindent \textsc{Proof:}
For each vertex $v$ of $\Gamma$ we have:
\begin{itemize}
 \item $f^{(j)}(v) \neq 0$ for all $j=1,\dots,d-1$, or
 \item $v$ has multiplicity at least $2$
\end{itemize}
(by the Gauss--Lucas theorem).
This means that among the $d$ roots of $f$, counting multiplicities,
at least $\gamma$ of them are not needed to
find a common root for each derivative.
If $\delta = 1$, some non-vertex has multiplicity $m$, so another $m-1$
roots are superfluous. Therefore, at most $d - \gamma - (m-1)$ roots are needed.
If $\delta = 0$, then the bound
reads $d - (\gamma - 1) - (m-1)$. In both cases, the upper bound for $\text{type}(f)$ follows.
The lower bound follows from an observation by
Draisma and Knopper~\cite[Proposition~6]{DD}.
\hfill $\blacksquare$\\

\noindent Refining to the level of scenarios, we find:

\begin{prop} \label{hightype}
Let $d > 2$ be an integer and let $s = (s_1, s_2, \dots, s_{d-1})$ be a scenario for degree $d$. If
\begin{enumerate}
 \item $\emph{type}(s) \in \{0,1,d-2\}$, or
 \item $\emph{type}(s) \leq d-3$, the first $d-2 - \emph{type}(s)$ entries of $s$ are zero,
   and among $s_{d-1-\emph{type}(s)}, \dots, s_{d-1}$ there is a zero or two consecutive entries that are equal,
\end{enumerate}
then there are no CA-polynomials $f \in \mathbb{C}[x]$ for which $\emph{scen}(f) = s$.
\end{prop}

\noindent \textsc{Proof:} The first part is an immediate corollary to Proposition~\ref{gonality}.
As for the second statement, suppose to the contrary that
$f$ is a CA-polynomial for which $\text{scen}(f) = s$, with
$t = \text{type}(s) \leq d-3$ and the first $d-2 - \text{type}(s)$ entries of $s$ equal to zero.
Let $a_0, \dots, a_t \in \mathbb{C}$ be as in (\ref{defscen}).
Then $a_0$ is a root with multiplicity at least $d-1-t$. Let $\Gamma$ be the convex
hull of the roots of $f$ and let $\gamma$ be its number of vertices. Using Proposition~\ref{gonality},
we conclude that $\gamma=2$ and that $a_0$ is a vertex. Then if another $0$ would appear in $s = \text{scen}(f)$,
by Gauss--Lucas we would conclude
that the multiplicity of $a_0$ is strictly bigger than $d-1-t$, which would contradict Proposition~\ref{gonality}.
On the other hand, if two consecutive entries would be equal, some high-order derivative of $f(x)$ would have a double root. But since
$\gamma = 2$, $f(x)$ is equivalent to a real-root polynomial, so Rolle's theorem would imply that this double root is actually
a root of $f(x)$ with multiplicity strictly bigger than $d - t$, again contradicting Proposition~\ref{gonality}.
 \hfill $\blacksquare$\\

\noindent \emph{Remark.} Let $s$ be as in the \'enonc\'e of Proposition~\ref{hightype}. Then one
cannot merely conclude (without using new arguments, that is) the stronger statement that there are no CA-polynomials $f \in \mathbb{C}[x]$ that
\emph{match} with $s$.\\

\noindent \textbf{(\ref{general}.3)}
As immediate corollaries to the lower bound $2 \leq \text{type}(f)$, we get the following three easy facts: if $f$ is
a CA-polynomial (over $\mathbb{C}$) of degree $d$, then
\begin{enumerate}
  \item $f^{(2)}(x)$ cannot be the $(d-2)$th power
  of a linear polynomial,
  \item $f$ cannot have a root of multiplicity at least $d-1$,
  \item $f$ has at least three distinct roots
\end{enumerate}
\noindent (note that these statements can be proved in various other ways, see e.g.\ \cite[Proposition 2.2]{Verh}).
In the next two propositions,
we will go a step further in directions 1 and 2. Later on (Proposition~\ref{3roots} and Theorem~\ref{4roots}), we will go two steps further
in direction 3.
\begin{prop}
If $f \in \mathbb{C}[x]$ is a CA-polynomial of degree $d$, then
$f^{(3)}(x)$ cannot be the $(d-3)$th power of a linear polynomial.
\end{prop}
\noindent \textsc{Proof:}
Suppose to the contrary that $f^{(3)}(x)$ is the $(d-3)$th power of a linear polynomial.
Thanks to Lemma \ref{change}, we may assume $f^{(3)}(x) = \frac{d!}{(d-3)!}x^{d-3}$.
Assume that $f^{(1)}(0)\not=0$, then $f$ has a root of multiplicity at least $2$ which is
different from $0$ and again by Lemma \ref{change}, we may
assume $f(1)=f^{(1)}(1)=0$.
Thus
\[f(x)=x^d-(d-1)x^2+(d-2)x;\ \ f^{(2)}(x)=(d-1)\left(dx^{d-2}-2\right).\]
Solving $f(x)=f^{(2)}(x)=0$, we get $x=\frac{d}{d+1}$ and $(\frac{d+1}{d})^{d-2}=\frac{d}{2}$.
We easily see that the function $\phi(t)=(t-2)\ln\frac{t+1}{t}-\ln\frac{t}{2}$ is strictly decreasing
for $t\ge 4$ and that $\phi(4)<0$. Thus the equality $\phi(d)=0$ is never reached for $d\ge 4$.
We conclude that we necessarily have  $f^{(1)}(0)=0$. Then, for some constant $c$,
$f^{(2)}(x)=d(d-1)x^{d-2}+2c$ and $f(x)=x^d+cx^2$. Solving $f(x)=f^{(2)}(x)=0$, we get that $c=0$.
\hfill $\blacksquare$

\begin{prop} \label{N-2}
Let $f \in \mathbb{C}[x]$ be a CA-polynomial of degree $d$, then $f$ cannot have a root of multiplicity at least $d-2$.
\end{prop}

\noindent \textsc{Proof:}
Suppose that $0$ is such a root. If $f^{(d-1)}(0)\not=0$, then we may assume that  $f(1)=f^{(d-1)}(1)=0$ and
\[f(x)=x^{d-2}(x^2-dx+d-1),\ \ f^{(d-2)}(x)=\frac{(d-1)!}{2}(dx^2-2dx+2).\]
Solving $f(x)=f^{(d-2)}(x)=0$, we get $x^2=2$ and $x=\frac{d+1}{d}$. Thus $(d+1)^2=2d^2$ which is impossible.
We conclude that we necessarily have
$f^{(d-1)}(0)=0$. Then, for some constant $c$, $f(x)=x^d+cx^{d-2}$ and
$f^{(d-2)}(x)=\frac{d!}{2}x^2+c$.
Solving $f(x)=f^{(d-2)}(x)=0$, we get $c=0$.
\hfill $\blacksquare$\\

\noindent We have chosen to present an elementary proof of Proposition~\ref{N-2}, though we also
can see it as a direct consequence of the forthcoming Proposition \ref{3roots}.\\

\noindent \textbf{(\ref{general}.4)} Let us recall some basic properties of the elementary symmetric polynomials.
  Let a polynomial $f$ and its derivatives be of the form
\[f^{(j)}(x)=\frac{d!}{(d-j)!} (x^{d-j}+{{d-j} \choose 1} a_1x^{d-j-1}+{{d-j} \choose 2} a_2x^{d-j-2}+\cdots+
a_{d-j})\]
(here by convention $f=f^{(0)}$).
Let $\sigma_m(j)$ be the sum of the $m$th powers of the roots of $f^{(j)}$, for $j=0,\cdots,d-1$. Then
Newton's formulas applied to each $f^{(j)}$ give the following relations
(see for example \cite{Pr} for more details on Newton formulas):

\begin{lem}\label{Relations}
\[\sum_{k=1}^r\sigma_k(j) {{d-j} \choose {r-k}}a_{r-k}=-r{{d-j} \choose r}a_r\]
for $0\le j\le d-1$, $1\le r\le  d-j$. (It is understood that $a_0=1$.)
\end{lem}

\noindent In particular, for $r=1$, we have that
\[\displaystyle \frac{\sigma_1(j)}{d-j}=\frac{\sigma_1(0)}{d}\]
 for  $j=0,\dots,d-1$, which means that the center of mass of the roots of the derivatives is fixed.
As obviously
\[\displaystyle \sigma_1(d-1)=\frac{\sigma_1(0)}{d}=-a_1\] is the only
root of $f^{(d-1)}$, we see that whenever $f$ is a
CA-polynomial over $\mathbb{C}$, the center of mass of its roots
$\displaystyle \frac{\sigma_1(0)}{d}$ is itself a root of $f$. As a direct consequence,
the number of distinct roots of a CA-polynomial cannot be two. Actually, we can say more:
if $f$ has more than two distinct roots, then at least one of them (the center of mass) has
to be in the interior of the convex hull of the roots. This fact also follows immediately from the Gauss--Lucas theorem,
and can be pushed further:

\begin{prop}\label{3roots}
Let $f \in \mathbb{C}[x]$ be a CA-polynomial. Then $f$ has at least two distinct roots in
the interior of the convex hull of the roots, when plotted in the complex plane. In particular, $f$ has at least four distinct roots.
\end{prop}

\noindent \textsc{Proof:}
Assume that $f$ has exactly one root, say $0$, in the interior. Let $\zeta$ be among the
roots of $f$ located on the boundary with maximal multiplicity $m$.
Then by Gauss--Lucas, $f^{(m)}(0)=f^{(m+1)}(0)=\cdots=f^{(d-1)}(0)=0$ which means that for $j=m,\ldots, d-1$:
\[\displaystyle f^{(j)}(x)=\frac{d!}{(d-j)!}x^{d-j}.\]
Taylor expansion gives
\[f(0)=\sum_{j=m}^d\frac{f^{(j)}(\zeta)}{j!}(-\zeta)^{j}=\zeta^d\sum_{j=m}^d (-1)^j{d \choose j}=\zeta^d (-1)^m{{d-1} \choose {m-1}}.\]
As $f(0)=0$, we get $\zeta=0$, which is a contradiction.
\hfill $\blacksquare$\\

\noindent Note that Proposition~\ref{3roots} can also be deduced directly from $2 \leq \text{type}(f)$.\\

\noindent \textbf{(\ref{general}.5)} We now prove the main result of this section:

\begin{thm}\label{4roots}
Let $f$ be a CA-polynomial over $\mathbb{C}$, then $f$ has at least five distinct roots.
\end{thm}

\noindent \textsc{Proof:}
 Assume that $f$ has four distinct roots. Then by the previous proposition, it has at least two distinct
roots in the interior of its Gauss--Lucas hull.
This implies that the four roots are on a line. By Lemma \ref{change}, we may assume that
this is the real line. We denote by $m$ the maximal multiplicity of the
roots of $f$. By Proposition~\ref{N-2}, we know that  $2\le m \le d-3$.
\begin{itemize}
\item First case: $m\le d-5$.
Again using Lemma \ref{change}, we may assume without loss of
generality that the roots of $f$ are as follows : $a<0<1<b$ and $f^{(d-1)}(0)=0$.
Then $a$ and $b$ cannot be zeros of $f^{(j)}$ for $d-5\le j\le d-1$. Moreover,
by Rolle's theorem, each zero of
$f^{(j)}$ is simple. Then we necessarily
have $f^{(d-2)}(1)=0,\ f^{(d-3)}(0)=0,\ f^{(d-4)}(1)=0, \ f^{(d-5)}(0)=0$.
Integrating five times the expression $f^{(d-1)}(x)=d!x$ and taking
into account these constraints, we
get $\displaystyle f^{(d-5)}(x)=\frac{d!}{5 !}x(x^2-5)^2$. But
this contradicts the fact that the roots are simple.

\item Second case: $m=d-4$. In view of
Lemma~\ref{change}, we arrange the roots as
follows : $a<0<b<1$ and we assume that $f^{(d-1)}(0)=0$.
Denote by $m_a,\ m_0,\ m_b,\ m_1$ their respective multiplicities.
Then again we must have $f^{(d-2)}(b)=0$, $f^{(d-3)}(0)=0$, $f^{(d-4)}(b)=0$.
Like in the first case, computing the last derivatives, we get
\[\begin{split}
&f^{(d-1)}(x)=d!x, \ \ \ 2! f^{(d-2)}(x)=d!(x^2-b^2),\\
&3! f^{(d-3)}(x)=d!x(x^2-3b^2),
 \ \ 4! f^{(d-4)}(x)=d!(x^2-5b^2)(x^2-b^2).
\end{split}
 \]
Obviously, as $f^{(d-4)}(b)=0$, we have $m_b\le d-5$.
From the Gauss--Lucas theorem, we deduce that $a<-\sqrt 5 b$. Now we apply
Lemma~\ref{Relations} with $j=0$, $r=1$ and with $j=0$, $r=3$ to obtain
\begin{equation}\label{equ}
m_a a+m_b b+m_1=m_a a^3+m_b b^3+m_1=0.
\end{equation}
We deduce that $m_a a(a^2-1)=-m_bb(b^2-1)$ and looking at the sign,
we see that $-a<1$. Then $m_a>-am_a=m_b b+m_1>m_1$ which implies that $m_a\ge 2$ and $m_1\le d-5$.
Now in the case where  $m_a=2, m_1=m_b=1$, equations~(\ref{equ}) give $a(a+1)^2=0$.
Thus this case cannot occur. We can readily deduce that $m_0\le d-5$.
The only possibility left is $m_a=m=d-4$.

From the relation
$-(d-4)a(1-a^2)=m_bb(1-b^2)$, we deduce that $\phi(-a)\le \phi(b)$ where we
put $\phi(t)=t(1-t^2)$. But $\phi$ is increasing on $[0,1/\sqrt 3]$ and we know
that $-a>b>0$. Thus we have $-a>1/\sqrt 3$. Now we get back to the
linear equation in (\ref{equ}):
\[d-4=m_b\frac{b}{-a}+m_1\frac{1}{-a}<\frac{m_b}{\sqrt 5}+m_1\sqrt 3<4.\]
Since the Casas-Alvero conjecture is true for $d \leq 7$, this is a contradiction.


\item Third case: $m=d-3$. We proceed as in the previous case.  We have
 \[
f^{(d-1)}(x)=d!x, \ \ \ 2! f^{(d-2)}(x)=d!(x^2-b^2), 3! f^{(d-3)}(x)=d!x(x^2-3b^2).
 \]
From Gauss--Lucas we deduce that $a<-\sqrt 3 b$. Again, we obtain  that $m_a\ge 2$. Thus we necessarily
have: $m_a=m$, $m_0=m_1=m_b=1$. The linear equation in (\ref{equ}) gives
\[d-3=\frac{b}{-a}+\frac{1}{-a}<\frac{1}{\sqrt 3}+\sqrt 3<3,\]
again a contradiction. \hfill $\blacksquare$
\end{itemize}

\section{Additional constraints for special degrees} \label{additional}

\noindent \textbf{(\ref{additional}.1)} We now turn
our attention to certain special instances of $d$, in each case involving a prime number $p$.
Inspired by Draisma and de Jong's take \cite{DD}, we use $p$-adic valuations.
Most of the proofs below have straightforward analogs in
the original reduction-mod-$p$ setting
of Graf von Bothmer et al. But at some points, the valuation language does seem slightly more powerful.
Our starting point is the existence of a map
\[ v_p  \, : \, \mathbb{C} \, \rightarrow \, \mathbb{Q} \cup \{ + \infty \} \]
satisfying
\begin{itemize}
  \item $v_p(a)=+\infty$ if and only if $a=0$,
  \item $v_p(ab)=v_p(a)+v_p(b)$ for all $a,b \in \mathbb{C}$,
  \item $v_p(a+b)\ge \min\{ v_p(a), v_p(b) \}$ for all $a,b \in \mathbb{C}$,
\end{itemize}
and extending the usual $p$-adic valuation on $\mathbb{Z}$ (i.e. if $n = p^r \cdot n'$ with $n'$ prime to $p$,
then $v_p(n) = r$). See e.g.\ \cite[Chapter 4, Theorem 1]{Rib}.
It is important to note that the last property implies
$v_p(a+b)= \min\{ v_p(a), v_p(b) \}$ if $v_p(a)\not=v_p(b)$. We will make a frequent use of this fact.\\

\noindent \textbf{(\ref{additional}.2)} The $p$-adic valuations of binomial
coefficients are well-understood. A formula due to Legendre \cite{Legendre}
states that for any $n \in \mathbb{Z}_{> 0}$ and any $j \in \{0, \dots, n\}$
one has
\[ v_p {n \choose j} = \frac{s_p(j) + s_p(n-j) - s_p(n)}{p-1}, \]
where $s_p(\cdot)$ denotes the sum of the $p$-adic digits.
Note that $s_p(j) + s_p(n-j) - s_p(n)$ is a measure for the number of carries
when adding $n-j$ to $j$ in base $p$. In particular,
\[ v_p {n \choose j} = 0 \quad \text{iff} \quad \text{there are no carries}.\]
It follows that:
\begin{lem} \label{binomvals}
 Let $n \in \mathbb{Z}_{>0}$ and $k \in \mathbb{Z}_{\geq 0}$. If $j \in \{0, 1, 2, \dots, np^k \}$ is
 not a multiple of $p^k$, then
 \[ v_p {np^k \choose j} > 0. \]
 If moreover $n = p^r + 1$ for some $r \in \mathbb{Z}_{\geq 0}$, it is sufficient to assume that
 $j \not \in \{ 0, p^k, (n-1)p^k, np^k \}$.
\end{lem}
\noindent \textsc{Proof:} According to Legendre's formula
\[ v_p {np^k \choose j} = \frac{s_p(j) + s_p(np^k-j) - s_p(np^k)}{p-1}. \]
Let $q$ and $\rho \neq 0$ be the quotient and remainder of $j$ when divided by $p^k$.
Then $s_p(np^k) = s_p(n)$, $s_p(j) = s_p(q) + s_p(\rho)$,
and
\[ s_p(np^k - j) = s_p((n - q - 1)p^k + (p^k-\rho)) \geq s_p(n-q) - 1 + 1, \]
from which
\[
 v_p {np^k \choose j} \geq v_p {n \choose q} + \frac{s_p(\rho)}{p-1} > 0.
\]
A similar argument proves the second statement. \hfill $\blacksquare$\\

\noindent \textbf{(\ref{additional}.3)} We use this to prove:
\begin{prop} \label{emptyscens}
 Let $n \in \mathbb{Z}_{>0}$ and $k \in \mathbb{Z}_{\geq 0}$ be integers,
 and let $f \in \mathbb{C}[x]$ be a CA-polynomial of degree $d = np^k$. Then
 \[ f, f^{(p^k)}, f^{(2p^k)}, \dots, f^{(d-p^k)} \]
 do not share a common root. If $n=p^r + 1$ for some integer $r \geq 0$,
 one even has that
 \[ f, f^{(p^k)}, f^{(d-p^k)} \]
 do not share a common root. As a consequence, if $s= (s_1, \dots, s_{d-1})$
 is a scenario for degree $d$ and $s_{p^k} = s_{2p^k} = \dots = s_{d-p^k}$
 (resp.\ $s_{p^k} = s_{d-p^k}$), then there are no CA-polynomials
 that match with $s$.
\end{prop}

\noindent \textsc{Proof:}
We only prove the first statement (the second assertion follows entirely similarly).
Suppose to the contrary that $f$ is a CA-polynomial such that $f, f^{(p^k)}, \dots, f^{(d-p^k)}$
do have a common root.
We may assume without loss of generality, using Lemma \ref{change}, that $f$ is of the form
\begin{equation} \label{standardform}
 f(x)=x^{d}+ {d \choose 1} a_1 x^{d-1} + {d \choose 2} a_2x^{d-2}+\cdots +  {d \choose d- 1} a_{d-1}x,
\end{equation}
that the assumed common root of $f, f^{(p^k)}, \dots, f^{(d-p^k)}$ is $0$, and that
\[ \min\{v_p(x_i) \, | \,  i=1,\dots,d\}=0,\]
where we have denoted by $x_1,x_2,\dots,x_d$ the zeros of $f$.

For $j=1,\dots,d-1$, we have:

\begin{equation}\label{derivee1}
\frac{j!}{d!}f^{(d-j)}(x)=x^{j}+{j\choose 1} a_1x^{j-1}
+ {j\choose 2} a_2 x^{j-2}+\cdots+{j\choose {j-1}}a_{j-1}x+a_j.
\end{equation}
Using equality (\ref{derivee1}) with $j=1,\cdots, d-1$,
each time plugging in
a common root of $f^{(d-j)}$ and $f$ (taking $0$ if $j$ is a multiple of $p^k$), one proves by induction on $j$  that
\begin{equation}\label{valuations6}
\left\{
\begin{array}{ll}
v_p(a_j)\ge 0 & \hbox{ for all\ }  \ j=1,\dots, d-1, \\
a_j = 0 & \hbox{ as soon as\ } \  p^k \mid j. \\
\end{array}
\right.
\end{equation}
Now let $x_j$ be such that $v_p(x_j)=0$. Then taking valuations of both sides of the equality
\[ x_j^d = - {d \choose 1} a_1 x_j^{d-1} - {d \choose 2} a_2x_j^{d-2} - \cdots  - {d \choose {d-2}} a_{d-2}x_j^2 - {d \choose d- 1} a_{d-1}x_j\]
yields a contradiction with (\ref{valuations6}) and Lemma~\ref{binomvals}. \hfill $\blacksquare$\\

\noindent Note that the cases $p^k$ and $2p^k$ tautologically follow from the above proposition. If $d = p^r + 1$,
it implies that the root of $f^{(d-1)}(x)$ must be a simple root of $f(x)$. If $p \geq 3$, this in turn
can be seen as a limit case of the following statement:

\begin{prop}
 If $d = p^r + 1$, then the root of $f^{(d-1)}(x)$ cannot be the mean of two distinct roots of $f(x)$.
\end{prop}

\noindent \textsc{Proof:} Using Lemma~\ref{change} we can assume that $f(x)$ is of the form
(\ref{standardform}) with $a_1 = 0$ (i.e.\ the root of $f^{(d-1)}(x)$ is $0$),
and that again all roots $x_1, \dots, x_d$ have non-negative valuation, with minimum $0$. Let
$x_j$ be such that $v_p(x_j)=0$. Then the equality
\[ da_{d-1}x_j = -x_j^d - {d \choose 2}a_2x_j^{d-2} - \cdots - {d \choose d-2} a_{d-2}x_j^2 \]
implies that $v_p(a_{d-1}) = 0$. Now let $w \in \mathbb{C}^\ast$ be such that
$f(w) = f(-w) = 0$. Then $0 = f(w)-f(-w)$ gives
\[ da_{d-1}w = - {d \choose 3}a_3w^{d-3} - {d \choose 5}a_5w^{d-5} - \cdots - {d \choose d-3}a_{d-3}w^3. \]
Taking valuations yields a contradiction. \hfill $\blacksquare$\\

\noindent The same argument can be used to show that the root of $f^{(d-1)}(x)$ cannot be the mean of two distinct
roots of $f^{(1)}(x)$.\\

\noindent \textbf{(\ref{additional}.4)} From now on, we focus on the special case $d = p+1$.
Using once again Lemma \ref{change}, we may assume that

\begin{equation}\label{p+1}
\left\{
\begin{array}{ll}
&f(x)=x^{d}+d a_1x^{d-1}+{d\choose 2} a_2x^{d-2}+\cdots+ {d\choose {d-2}} a_{d-2}x^2,\\
& \\
&\min\{v_p(x_j) \, | \, j=1,\dots,d\}=0, \end{array}
\right.
\end{equation}
where we have denoted by $x_1,\dots, x_{d-3}, x_{d-2}=x_{d-1}=0, x_d=-a_1$ the roots of $f$.
For $j=1,\dots,d-2$, we then again have that expression (\ref{derivee1}) holds.
Observe that $v_p(a_1)\ge 0$ because $-a_1$ is one of the roots of $f$.
As before, using equality (\ref{derivee1}) with $j=2,\dots, d-2$,
each time plugging in a common root of $f^{(d-j)}$ and $f$, we prove by induction on $j$  that
\begin{equation}\label{valuations2}
v_p(a_j)\ge 0 \ \ \ \hbox{ for all\ }  \ j=1,\dots, d-2.
\end{equation}
Let $x_j$ be such that $v_p(x_j)=0$.  The equality
\[-d a_1 x_j^{d-1}= x_j^d+{d\choose 2} a_2x_j^{d-2}+\cdots+ {d\choose {d-2}}a_{d-2} x_j^2\]
shows that $v_p(a_1)=0$. Therefore, we may assume without loss of generality that $a_1=-1$. Then we can write
$f(x)=(x-1)g(x)$ where
\[
\begin{split}
g(x)&=x^{d-1}-(d-1)x^{d-2}+ \left( {d\choose 2} a_2-(d-1) \right)x^{d-3}+\\
& \left( {d\choose 3} a_3+{d\choose 2}a_2-(d-1) \right)x^{d-4}+\cdots +\\
& \left( {d\choose d-3}a_{d-3}+\cdots+{d\choose {2}} a_2-(d-1) \right)x^2.
\end{split}\]
In view of (\ref{valuations2}) and Lemma \ref{binomvals},
all roots of $g$ have strictly positive valuations  (actually greater than $1/(d-3)$).
As a consequence, we see that $1$ is a simple root of $f$ (a fact already implied by Proposition~\ref{emptyscens})
and that $v_p(x_j)>0$ for $j=1,\dots,d-3$.
Now whenever $f^{(d-j)}(1)\not=0$, the
Casas-Alvero property implies that $f^{(d-j)}(x_j)=0$
with $v_p(x_j)>0$ and from equality (\ref{derivee1}) we
get $v_p(a_j)>0$.
But as
\[ f(1)=1-d+{d\choose 2}a_2+\cdots+{d\choose {d-2}} a_{d-2}=0,\]
there is at least one index $2\le j\le d-2$ such that $v_p(a_j)=0$.
In other words, at least one of the derivatives $f^{(d-j)}(1)=0$.
If we put this together with Proposition~\ref{N-2} and the observations
following Lemma~\ref{Relations}, we get:
\begin{lem}\label{1fois}
Let $f$ be a CA-polynomial over $\mathbb{C}$ of degree $d=p+1$, where $p$ is
prime. Let $c$ be the center of mass of the roots of $f$. Then the  following conditions are satisfied:
\begin{itemize}
\item $f^{(1)}(c) \neq 0, f^{(d-1)}(c)=0$,
\item $f^{(j)}(c)\not=0$ for at least one $j\in \{2,\dots,d-2\}$,
\item $f^{(j)}(c)=0$ for at least  one $j\in \{2,\dots,d-2\}$.\\
\end{itemize}
\end{lem}

\noindent \textbf{(\ref{additional}.5)} Let us now go further into the
investigation of the orders of the derivatives having the center of mass as
a root, thereby proving Theorem~\ref{DeltaThm}. We may again assume that $f$ is of the form (\ref{p+1}) and that $a_1=-1$.
We will use the notation $x\equiv y$ if $v_p(x-y)>0$.
In view of Lemma \ref{1fois}, let $j_1<j_2<\cdots<j_m\ $ be the indices between $2$ and $d-2$ such that
$f^{(d-j_i)}(1)=0$ for $i=1,\dots,m$.
As observed previously, for all $j\in\{2,\cdots,d-2\}$, we have $v_p(a_j)\ge 0$.
Moreover, if $j\notin \{j_1,\cdots, j_m\}$ then $a_j\equiv 0$.
From equality (\ref{derivee1}) with $x=1$ and $j=j_1, j_2,\dots, j_m$, we get
\begin{equation}\label{S}
\left\{
\begin{array}{ll}
1-j_1+a_{j_1} & \equiv 0\\
1-j_2+\binom{j_2}{j_1}a_{j_1}+a_{j_2} & \equiv 0\\
\ \ \ \ \  \vdots\\
1-j_m+\binom{j_m}{j_1} a_{j_1}+\binom{j_m}{j_{2}}a_{j_2} +\cdots+ a_{j_m} & \equiv 0\\
\end{array}
\right.
\end{equation}
Now, using that  $\frac{f(1)}{p}=0$ and that $v_p\binom{d}{j}\ge 1$ for $j=2,\dots, d-2$, we obtain
\begin{equation}\label{f(1)=0}
-1+\frac{\binom{d}{j_1}}{p} a_{j_1}+\cdots+\frac{\binom{d}{j_m}}{p}a_{j_m} \equiv 0.
\end{equation}
Observe that
 for all $2\le j\le d-2$ we have:
\[\begin{split}
\frac{{d\choose j}}{p}&=\frac{d(d-2)(d-3)\cdots(d-(j-1))}{j!}\\
&=\frac{(p+1)(p-1)(p-2)\cdots(p-(j-2))}{j!}\\
&=\frac{1}{j!}(p^{j-1}+\alpha_{j-2}p^{j-2}+\cdots+
\alpha_1p)+\frac{(-1)^{j-2}(j-2)!}{j!}\end{split}\]
where $\alpha_1,\dots, \alpha_{j-2}$ are integers.
Therefore: \[ \frac{\binom{d}{j}}{p}\equiv \frac{(-1)^j}{j(j-1)}.\]
Putting equations (\ref{S}) and (\ref{f(1)=0}) together and putting $\tilde a_{j_i}=\frac{a_{j_i}}{j_i(j_i-1)}$, we obtain:
\begin{equation}\label{systeme}
\left\{
\begin{array}{ll}
-1+j_1\tilde a_{j_1} & \equiv 0\\
-1+\binom{j_2-2}{j_1-2}j_2 \tilde a_{j_1}+j_2\tilde a_{j_2} & \equiv 0\\
\ \ \ \ \  \vdots\\
-1+\binom{j_m-2}{j_1-2}j_m \tilde a_{j_1}+\binom{j_m-2}{j_{2}-2}j_m \tilde a_{j_2} +\cdots+ j_m \tilde a_{j_m} & \equiv 0\\
-1+(-1)^{j_1} \tilde a_{j_1}+(-1)^{j_2}\tilde a_{j_2}+\cdots+(-1)^{j_m}\tilde a_{_m} & \equiv 0.
\end{array}
\right.
\end{equation}
With $\Delta_f$ as in the \'enonc\'e
of Theorem~\ref{DeltaThm}, we
see that necessarily $\det \Delta_f \equiv 0$: otherwise
inverting (\ref{systeme}) we would get that $1\equiv 0$. To conclude the proof of Theorem~\ref{DeltaThm} we show:
\begin{lem}\label{2fois}
Let $f \in \mathbb{C}[x]$ be a CA-polynomial of degree $d=p+1$ and let $c$ be the center of mass of its roots.
Then there are at least two indices $2\le j_1<j_2\le d-2$ such that
$f^{(j_1)}(c)=f^{(j_2)}(c)=0$.
\end{lem}
\noindent \textsc{Proof:}
If not, in virtue of Lemma \ref{1fois}, there exists a unique index $2\le j\le d-2$ such that $f^{(d-j)}(c)=0$.
We can assume without loss of generality that $f$ is of the form (\ref{p+1}) with $a_1=-1$ and apply the above.
Then $m=1$ and
\begin{equation}
\Delta_f=\left[
\begin{array}{l l }
-1& j \\
-1 & (-1)^{j}
\end{array}
\right]=j-(-1)^j.
\end{equation}
Observe that $1\le j-(-1)^j\le j+1\le d-2$ for $j\in{2,\dots, d-3}$.
Besides, $d-2-(-1)^{d-2}=d-3$ because $d$ is even (indeed, $p \neq 2$ since
the Casas-Alvero conjecture is true for degree $3$).
Thus there is no way for $p$ to divide $\det \Delta_f$. \hfill $\blacksquare$\\

\noindent \textbf{(\ref{additional}.6)} Theorem~\ref{DeltaThm} implies that every
CA-polynomial of degree
$d=p+1$ matches with a scenario $s=(s_1, \dots, s_{d-1})$
for which $s_{d-1} \neq 0$ and the index set
\[ \text{ind}(s) = \{ \, j \ | \ 2 \leq j \leq d-2 \text{ and } s_{d-j} = s_{d-1} \, \} \]
satisfies the according determinant condition. We remark however that
this does not necessarily imply that \emph{the} scenario
of a CA-polynomial satisfies these conditions. Indeed, imagine
a CA-polynomial $f \in \mathbb{C}[x]$ of degree $12$ for which
\[\text{scen}(f) = s = (0,1,2,3,4,2,5,6,4,7,4), \] i.e.\ there exist
$a_1, \dots, a_7 \in \mathbb{C}$ such that $f(a_{s_j}) = f^{(j)}(a_{s_j}) = 0$
for $j=1, \dots, d-1$. Then $\text{ind}(s) = \{3,7\}$ does not satisfy the determinant condition.
However, it might a priori be that $f^{(6)}(x)$ has both $a_2$ and $a_4$ as a root.
Then $f(x)$ also matches with
the scenario $(0,1,2,3,4,4,5,6,4,7,4) \neq \text{scen}(f)$. Here, the index set reads $\{3,6,7\}$, for
which the determinant condition \emph{is} satisfied.\\

\noindent \textbf{(\ref{additional}.7)} We end our study of the degree $p+1$ case
with the following observation.
\begin{prop}\label{rational}
Let $p$ be a prime number. Then there is no CA-polynomial of degree $d=p+1$ all of whose roots are rational.
\end{prop}

\noindent \textsc{Proof:}
Using the notations and the results found in the proof of Lemma \ref{1fois}, we may assume that $f$ is of the form
\[
\begin{split}
f(x)= & x^{d}-d x^{d-1}+{d\choose 2}x^{d-2}\\
    &+\cdots+(-1)^{k-1}{d\choose {k-1}}x^{d-k+1}+{d\choose k}a_kx^{d-k}+\cdots+{d\choose d-2} a_{d-2}x^2,
    \end{split}\]
    with
 $v_p(x_j)\ge 1$ for $j=1,\dots, d-3$.
Here, we have denoted by $k$ the smallest index between $2$ and $d-2$ such that $f^{(d-k)}(1)\not=0$
(we know from Lemma~\ref{1fois} that such a $k$ exists).
 We introduce the notation
 \[S_m=\sum_{j=1}^{d-3} x_j^m.\]
 Then we have: $v_p(S_1)=v_p(d-1)=1$, and $v_p(S_j)\ge 2$ for
$j=2,\dots,d-2$. Using Newton's formulas (see Lemma \ref{Relations} applied to $j=0$), we obtain
\[
\begin{split}
 -k{d\choose k}a_k&=\sum_{j=0}^{k-1}(-1)^j(1+S_{k-j}){d\choose j}\\   &=\sum_{j=0}^{k-1}(-1)^j{d\choose j}+\sum_{j=0}^{k-1}(-1)^jS_{k-j}{d\choose j}\\
& =(-1)^{k-1}{{d-1}\choose {k-1}}+\sum_{j=0}^{k-1}(-1)^jS_{k-j}{d\choose j}.\\
 \end{split}\]
 Note that $v_p({d\choose k}a_k)>1$ which will lead to a contradiction:
\begin{itemize}
 \item If $k=2$, then the last equality becomes
\[-2{d\choose 2} a_2=-(d-1)+S_2-dS_1=-(d-1)+S_2-d(d-1)=-(d+1)(d-1)+S_2.\]
The valuation of the right-hand term is $1$.
 \item If $3\le k\le d-2$, then the right-hand term is
 \[(-1)^{k-1}{{d-1}\choose {k-1}}+\sum_{j=0}^{k-2}(-1)^jS_{k-j}{d\choose j}+(-1)^{k-1}S_1{d\choose {k-1}}.\]
 But $v_p(S_{k-j})\ge 2$ for $j=0,\dots,k-2$, and $v_p(S_1{d\choose {k-1}})=2$, so the valuation of the right-hand term
 is $v_p({{d-1}\choose {k-1}})=1$. \hfill $\blacksquare$
\end{itemize}

\noindent Remark that the proof of Proposition~\ref{rational} in fact
implies that there are no CA-polynomials of degree $p+1$ all of whose
roots are contained in a number field in which $p$ does not ramify. Indeed,
this ensures that the valuations of the $x_j$ are integers, hence we
can still conclude that $v_p(x_j) \geq 1$.

\section{Algebraic varieties of counterexamples} \label{geometric}

\noindent \textbf{(\ref{geometric}.1)}
Let $k$ be an algebraically closed field and let $d>0$ be an integer. The set of equivalence classes
(in the sense of Lemma~\ref{change}) of CA-polynomials of
degree $d$ will be denoted by $\text{CA}_k(d)$.\\

\noindent \textbf{(\ref{geometric}.2)} We have a surjective map
\[ \Phi_k(d,d-2) : V_k(d,d-2) \rightarrow \text{CA}_k(d) : (p_1, \dots, p_{d-2}) \mapsto x^2(x-p_1)\cdots (x-p_{d-2}), \]
where $V_k(d,d-2) \subset \mathbb{P}_k^{d-3}$ is the projective variety defined by the ideal
\[ I_k(d,d-2) = \left( \, \left. \text{Res}_x(F,F_H^{(j)}) \, \right| \, j = 2, \dots, d-1 \, \right) \]
with $F = x^2(x - P_1) \dots (x - P_{d-2}) \in k[P_1,\dots,P_{d-2}][x]$. Therefore,
in order to prove that no CA-polynomials exist in degree $d$, it suffices to show that
$V_k(d,d-2) = \emptyset$. Note
that $V_k(d,d-2)$ is invariant under coordinate permutations, so it is sufficient to show that $V_k(d,d-2)$ does not contain
any points of the form $(p_1, \dots, p_{d-3},1)$. Setting $P_{d-2} = 1$ in $I_k(d,d-2)$, we obtain an ideal
of $k[P_1,\dots,P_{d-3}]$ that is equal to the unit ideal if and only if $V_k(d,d-2) = \emptyset$. This can be checked using a finite Gr\"obner basis
computation, which is exactly the approach of \cite{DTGV}.\\

\noindent \textbf{(\ref{geometric}.3)} Somehow dually, we also have a surjective map
\[ \Phi_k(d,0) : V_k(d,0) \rightarrow \text{CA}_k(d) : \qquad \qquad \qquad \qquad \qquad \]
\[ \qquad \qquad \qquad (a_1, \dots, a_{d-2}) \mapsto x^2(x^{d-2} + a_1x^{d-3} + \dots + a_{d-2}), \]
where now $V_k(d,0) \subset \mathbb{P}_k(d-2; d-1; \dots; 2; 1)$ is the weighted projective variety defined by the ideal
\[ I_k(d,0) = \left( \, \left. \text{Res}_x(F,F_H^{(j)}) \, \right| \, j = 2, \dots, d-1 \, \right) \]
with $F = x^2(x^{d-2} + A_1x^{d-3} + \dots + A_{d-2}) \in k[A_1,\dots,A_{d-2}][x]$. Again, in order to show that
no Casas-Alvero polynomials can exist in degree $d$, it is sufficient to prove that $V_k(d,0) = \emptyset$. This was used
in the theoretical approach of \cite{GLSW}.\\

\noindent \textbf{(\ref{geometric}.4)} We will make use of a hybrid version of the above maps. Namely, for
each $t \in \{ 0, \dots, d-2 \}$ we have a surjective map
\[ \Phi_k(d,t) : V_k(d,t) \rightarrow \text{CA}_k(d) : \qquad \qquad \qquad \qquad \qquad \qquad \qquad \]
\[ (p_1, \dots, p_t,a_1,\dots,a_{d-2-t}) \mapsto x^2(x-p_1)\cdots (x-p_t)(x^{d-2-t} + a_1x^{d-3-t} + \dots + a_{d-2-t}), \]
where $V_k(d,t) \subset \mathbb{P}_k(1;\dots;1;d-2-t;d-3-t;\dots;2;1)$ is the weighted projective variety defined by the ideal
\[ I_k(d,t) = \left( \, \left. \text{Res}_x(F,F_H^{(j)}) \, \right| \, j = 2, \dots, d-1 \, \right) \]
with
\[ F = x^2(x-P_1)\cdots (x-P_t)(x^{d-2-t} + A_1x^{d-3-t} + \dots + A_{d-2-t})\]
in $k[P_1, \dots, P_t,A_1,\dots,A_{d-2-t}][x]$.
Once more it is sufficient to show that $V_k(d,t) = \emptyset$ (for any value of $t$) in order
to prove that no Casas-Alvero polynomials of degree $d$ exist over $k$.\\

\noindent \textbf{(\ref{geometric}.5)} Now to each scenario $s$ for degree $d$ of type $t$, we associate the variety
\[ V_k(s) \subset V_k(d,t) \]
defined by the ideal
\[ I_k(s) = \left( \, \left. F_H^{(j)}(P_{s_j}) \, \right| \, j = 2, \dots, d-1 \, \right)
\subset k[P_1, \dots, P_t,A_1,\dots,A_{d-2-t}] \]
where
\[ F = x^2(x-P_1)\cdots (x-P_t)(x^{d-2-t} + A_1x^{d-3-t} + \dots + A_{d-2-t})\]
and $P_0 = 0$.
Then it is clear that $V_k(s)$ parameterizes the CA-polynomials
that match with $s$. Recall that every CA-polynomial matches with at least
one scenario (e.g., its own scenario $\text{scen}(f)$). Thus, if one wants to show that no CA-polynomials
of degree $d$ exist over $k$, it suffices to show that $V_k(s) = \emptyset$
for each scenario $s$ for degree $d$. 
This is essentially the `primary decomposition' that was mentioned in \cite[Section~?]{GLSW},
but in Section~\ref{computational} below 
we will see that there is a significant amount of computational gain to be expected
from viewing the set of CA-polynomials that match with $s$ as a 
subvariety of $V_k(d,t)$ rather than $V_k(d,d-2)$.
Moreover, if $k = \mathbb{C}$, in view of the theoretical results obtained
in Sections~\ref{general} and~\ref{additional}, it is actually sufficient
to check whether $V_\mathbb{C}(s) = \emptyset$ for a restricted set of scenarios. We will elaborate the details of
this for $d=12$ in Section~\ref{degree12}.

\section{Revisiting the computational approach} \label{computational}

\noindent \textbf{(\ref{computational}.1)} We now
describe the basic version of our algorithm, discarding
the theoretical results of Sections~\ref{general} and~\ref{additional}. The input is a field characteristic $p$ (either $0$ or
a prime number) along with an integer $d > 2$. The output is \verb"yes" or \verb"no", depending on whether
Casas-Alvero polynomials exist in degree $d$ and characteristic $p$ or not.\\

\noindent \emph{Step 1.} Create a list $L$ (of length $d-1$) of lists, such
that $L[t]$ contains all scenarios for type $t$ (for $t=0, \dots, d-2$). This can be
done easily using $d-2$ nested for-loops. Let $k$ be the field of rational numbers if $p=0$, and let
$k$ be the field with $p$ elements otherwise. Set \verb"answer := no".\\

\noindent \emph{Step 2.}  For $t$ going from $1$ to $d-2$ do:
\begin{itemize}
\item[-] Initiate the following variables/structures:
\begin{itemize}
 \item[*] $R = k[P_1, \dots, P_{t-1}, A_1, \dots, A_{d-2-t}]$
 \item[*] $S = R[x]$
 \item[*] $P_0 = 0$ and $P_t = 1$
 \item[*] $F(x) = x^2(x-P_1)\cdots(x-P_t)(x^{d-2-t} + A_1x^{d-3-t} + \dots + A_{d-2-t})$
 \item[*] $\prec$ $=$ a monomial ordering that first eliminates $A_1, \dots, A_{d-2-t}$ and that behaves like
\verb"grevlex" on the remaining variables $P_1, \dots, P_{t-1}$
\end{itemize}
\item[-] For $s$ in $L[t]$ do:
\begin{itemize}
\item[*] Let $I_k^\text{aff}(s) \subset R$ be the ideal generated by $F_H^{(j)}(P_{s_j})$ for $j=2, \dots, d-1$.
Check whether or not $I_k^\text{aff}(s) = R$ by checking if the reduced Gr\"obner basis (w.r.t.\ $\prec$)
of $I_k^\text{aff}(s)$ equals $\{ 1 \}$. If it does not, set \verb"answer := yes" and quit the loops.
\end{itemize}
\end{itemize}

\noindent \emph{Step 3.} Output \verb"answer".\\

\noindent \textbf{(\ref{computational}.2)} Modulo a base change to the algebraic
closure of $k$, $I_k^\text{aff}(s)$ is obtained from $I_k(s)$ (as described in \textbf{(\ref{geometric}.5)})
by setting $P_t = 1$, so it only describes an affine part of $V_k(s)$. However, it suffices to verify
that this affine part is empty.
Indeed, the type of a CA-polynomial corresponding to a point $(p_1, \dots, p_t, a_1, \dots, a_{d-2-t}) \in V_k(s)$
with $p_t = 0$ is strictly smaller than $t$, so we would have encountered it already.\\

\noindent \textbf{(\ref{computational}.3)} The variables $A_1, \dots, A_{d-2-t}$ appear linearly in the
defining polynomials $F_H^{(j)}(P_{s_j})$. Therefore, they can be eliminated easily. (In fact, the corresponding
linear system is in echelon form, so the $A_i$'s can be eliminated bottom-up by hand.) The lower the type,
the more variables can be eliminated and the easier the Gr\"obner basis computation becomes
(in the extreme case $t=1$ one obtains a linear system in $d-3$ variables).
This is the main reason for our usage of the hybrid varieties $V_k(d,t)$.\\

\noindent \textbf{(\ref{computational}.4)} It is theoretically possible to avoid Gr\"obner basis computations and use linear algebra instead.
Indeed, $I_k^\text{aff}(s) = R$ is equivalent to the solvability of
\begin{equation} \label{nullstellen}
 1 = g_1 \cdot F_H^{(2)}(P_{s_2}) \, + \, \dots \, + \, g_{d-2} \cdot F_H^{(d-1)}(P_{s_{d-1}})
\end{equation}
in terms of polynomials $g_i \in R$. If such polynomials exist, by the effective Nullstellensatz
they can be chosen such that their degree is
bounded by $d^d$ (e.g., see \cite{Kollar}). So in principle, one could use indetermined coefficients to translate
the solvability of (\ref{nullstellen}) to the solvability of some linear system of equations. But this system is so huge
that no gain is to be expected (although maybe this deserves a deeper analysis).\\

\noindent \textbf{(\ref{computational}.5)} One can speed up the algorithm slightly by noting the following. If $s_2 = 0$, then
the first defining polynomial is
\[ F_H^{(2)}(0) = (-1)^t \cdot P_1 \cdots P_{t-1} \cdot A_{d-2-t} \]
But Casas-Alvero polynomials corresponding to $P_1 \cdots P_{t-1} = 0$ are of strictly lower type
than $t$, so they would have been encountered already. Therefore, our defining polynomial can be replaced by $A_{d-2-t}$.
If in addition $s_3 = 0$, then similarly the second defining polynomial can be replaced by $A_{d-3-t}$, and so on. Suppose that
the first nonzero entry of $s$ appears at position $j$. Then after substituting $A_{d-2-t} = \dots = A_{d-j+1-t} = 0$
(no substitutions if $j=2$), one finds that
\[ F_H^{(j)}(P_{s_j}) = F_H^{(j)}(P_1) \]
is a multiple of $P_1$. For the same reason, this factor can be removed.\\

\noindent \textbf{(\ref{computational}.6)}
The above algorithm can be used straightforwardly to find all bad primes for a given degree $d$ (given that we know
that the Casas-Alvero conjecture is true in degree $d$):
\begin{enumerate}
 \item Initialize a set of candidate bad primes $C = \{ \, \}$.
 \item First run the basic algorithm with $p=0$, but instead of just checking whether the reduced Gr\"obner
basis of $I_\mathbb{Q}^\text{aff}(s)$ equals $\{1\}$, compute polynomials $g_1, \dots, g_{d-2} \in R$
for which (\ref{nullstellen}) holds. Then add every prime factor appearing in the denominators of the $g_j$ to $C$.
 \item Now if a prime $p$ is not in $C$, it cannot be a bad prime because each of the expansions (\ref{nullstellen})
 can be reduced mod $p$. To find which candidate bad primes are actually bad primes, we run the basic algorithm for each $p \in C$.
\end{enumerate}
An implementation of this method can be found in \verb"CAbadprimes.m".\\

\noindent \textbf{(\ref{computational}.7)}  The hardest part is step $2$, because of the computing in characteristic $0$.
Note that it is possible to give an upper bound for the elements of $C$ purely in terms of $d$, so that
step $2$ could in principle be avoided. Indeed, see the discussion following (\ref{nullstellen}) -- the denominators of the solutions
of the linear system can be bounded using Cramer's rule. But the bound one obtains is too large to be of any practical use.\\

\noindent \textbf{(\ref{computational}.8)}  We have executed the algorithm for $d=5$, $d=6$ and $d=7$. In case of $d=5$,
the total time needed was less than $0.03$ seconds. For $d=6$, the computer needed less than $3$ seconds.
A naive run of the algorithm for $d=7$ is not expected to end in a
reasonable amount of time, because the denominators become very hard to factor.
But by using several monomial orders and computing greatest common divisors, one can make the case $d=7$ feasible in Magma
(apart from the factorization of one composite $119$-digit number, for which we used the \verb"CADO-NFS" package \cite{CADO}).
The file \verb"CAbadprimes7test.m" contains Magma code proving the correctness of our output.
The case $d=8$ lies out of reach. Of course, exhaustive lists of bad primes
for increasing degrees become less and less interesting. But it would be good to have an idea on the growth of the largest
bad prime, or on the number of bad primes. Such lists can also be helpful in detecting patterns
(we could not observe any). By just repeating our basic algorithm for increasing values of $p$,
it is feasible to find the smallest non-bad prime (that does not divide $d$), for $d$ up to $10$.
We have put the outcomes in Table~\ref{nonbadprimes}.\\

\begin{table}
\begin{center}
\begin{tabular}{|l||l|l|l|l|l|l|l|l|l|}
\hline
$d$ & 2 & 3 & 4 & 5 & 6 & 7 & 8 & 9 & 10 \\
\hline
$p$ & - & - & 11 & 13 & 17 & 127 & 419 & 941 & 3803 \\
\hline
\end{tabular}
\caption{The smallest non-bad prime $p$ that does not divide $d$}
\label{nonbadprimes}
\end{center}
\end{table}

\section{The Casas-Alvero conjecture in degree $12$} \label{degree12}

\noindent \textbf{(\ref{degree12}.1)}
Naively applying the basic algorithm to $d=12$ and characteristic $p=0$ is unrealistic.
Two observations lead to a crucial speed-up:
\begin{itemize}
  \item as remarked in \textbf{(\ref{geometric}.5)}, in view of the theoretical results obtained in Sections~\ref{general} and~\ref{additional},
  it suffices to
        show that $V_\mathbb{C}(s) = \emptyset$ for a restricted set of scenarios $s$,
  \item for each such $s$, it actually suffices to show that $V_{\overline{\mathbb{F}}_p}(s) = \emptyset$ for
        a single prime $p$, because the varieties are projective and take equations over $\mathbb{Z}$.\\
\end{itemize}

\noindent \textbf{(\ref{degree12}.2)}
As for the first speed-up, by Theorem~\ref{DeltaThm} and Proposition~\ref{emptyscens}
it suffices to prove that $V_\mathbb{C}(s) = \emptyset$ for all scenarios $s = (s_1, \dots, s_{11})$
for which
\begin{itemize}
  \item $s_1 = 0 \neq s_{11}$,
  \item $s_3 \neq s_9$,
  \item $s_4 \neq s_8$,
  \item $\text{ind}(s)$ satisfies the determinant condition mentioned in the \'enonc\'e of
  Theorem~\ref{DeltaThm}
\end{itemize}
(we omit the contribution of Proposition~\ref{hightype} to this discussion, because
the arguments involved are rather subtle, whereas the computational gain is limited).
Let $L_{\text{res}}$ be obtained from $L$ (as introduced in \textbf{(\ref{computational}.1)})
by restricting to these scenarios. Then $\text{L}_{\text{res}}$ contains
\[ 0,6,718,5210,8918,5404,1352,141,5,0,0 \]
scenarios of type $0, \dots, 10$, respectively (this is less than was mentioned in (\ref{restrictedlist}), where
only the determinant condition was taken into account).
However, for the algorithm to work rigorously, the list $L_{\text{res}}$ should be slightly enlarged again, so that it becomes
closed under taking \emph{descendants}, in the following sense.
\begin{defn}
  Let $d > 0$ be an integer and let $s = (s_1, \dots, s_{d-1})$
  be a scenario for degree $d$. Let $t = \text{type}(s)$. Then we
  say that $s'=(s'_1, \dots, s'_{d-1})$ is a \emph{descendant} of $s$ if there exists a $1 \leq j \leq t$
  such that for all $i = 1, \dots, d-1$
  \begin{itemize}
    \item $s'_i = s_i$ if $s_i < j$,
    \item $s'_i = 0$ if $s_i = j$,
    \item $s'_i = s_i-1$ if $s_i > j$.
  \end{itemize}
\end{defn}
\noindent This ensures that working in the affine subvariety $P_t = 1$ (see
\textbf{(\ref{computational}.2)}) and speeding up the algorithm (as in
\textbf{(\ref{computational}.5)}) are still justified. Note that if $s'$ is a descendant
of $s$, then $\text{type}(s') = \text{type}(s) - 1$. By closing $L_{\text{res}}$ under taking
descendants, one obtains a list $L_{\text{res}}^{\text{cl}}$ containing
\[ 1, 279, 3892, 12073, 13661, 6685, 1491, 146, 5, 0, 0 \]
scenarios of type $0, \dots, 10$, respectively. This may seem a big increase, but note that
scenarios of low type can be eliminated very easily.\\


\noindent \textbf{(\ref{degree12}.3)}
As for the second speed-up, based on the experimentally observed distribution of bad primes
in degrees $d \leq 7$, any prime $p$ which is `not too small' is most likely to work.
If nevertheless the computation breaks down and a \verb"yes" is printed, one can redo the computation using a different value of $p$.
(In principle, it is possible to give a lower
bound on $p$ so that it is guaranteed to work, but this bound is much too large
to be of any practical use -- recall from Theorem~\ref{badprimesoverview} that the largest
bad prime for $d=7$ had already $135$ decimal digits). Our first try was $p = 10^7 + 17$ and immediately worked.
It is convenient to use the same $p$ for all scenarios listed in $L_{\text{res}}^{\text{cl}}$.
At least, if a scenario $s$ is treated modulo some $p$, then all of its subsequent descendants should
be treated modulo the same $p$. Indeed, this enables us to conclude that the \emph{projective} variety
$V_{\overline{\mathbb{F}}_p}(s)$ is empty, and hence that $V_\mathbb{C}(s) = \emptyset$.\\

\noindent \textbf{(\ref{degree12}.4)} Magma code implementing the above method can be found
in the file \verb"CAdeg12.m". We have executed the algorithm and the outcome was affirmative (i.e.\ the Casas-Alvero
conjecture is true in degree $12$, thereby proving Theorem~\ref{deg12thm}). Approximate time and memory requirements
can be found in Table~\ref{timings}.\\
\begin{table}
\begin{center}
\begin{tabular}{|l|l|l|l|}
\hline
type & \# scenarios & time & memory \\
\hline
$1$ & $279$ & $0.1$ secs & $\ll 0.1$ GB \\
$2$ & $3892$ & $43$ secs & $\ll 0.1$ GB \\
$3$ & $12073$ & $2$ mins & $< 0.1$ GB \\
$4$ & $13661$ & $40$ mins & $0.1$ GB \\
$5$ & $6685$ & $20$ hours & $0.2$ GB \\
$6$ & $1491$ & $2$ weeks & $1.3$ GB \\
$7$ & $146$ & $16$ weeks & $10$ GB \\
$8$ & $5$ & $15$ weeks & $90$ GB \\
\hline
\end{tabular}
\caption{Approximate time and memory requirements for settling $d=12$,
as if the algorithm were executed on a single core. In practice, types $6$ and $7$
were spread among multiple cores. In case of type $8$, this was not possible due to memory limitations.
}
\label{timings}
\end{center}
\end{table}

\noindent \textbf{(\ref{degree12}.5)} The computation fills in the smallest open entry in the
list of degrees for which the Casas-Alvero conjecture is known to hold. Up to our knowledge, the list of
degrees $d \leq 100$ for which the conjecture is still open is
\begin{center}
\footnotesize{$20,24,28,30,35,36,40,42,45,48,55,56,60,63,66,70,72,77,78,80,84,88,90,91,98,99,100$}.
\end{center}
Our algorithm can in principle be generalized to higher degrees (note in particular that the
two next open cases $d=20$ and $d=24$ are also of the form $p+1$). But without
new theoretical ingredients, an implementation of this is expected to demand astronomical
amounts of time and memory.

\small

\noindent \textsc{Departement Wiskunde, KU Leuven}\\
\noindent \textsc{Celestijnenlaan 200B, 3001 Leuven (Heverlee), Belgium}\\
\noindent \emph{E-mail address:} \verb"wouter.castryck@wis.kuleuven.be"\\

\noindent \textsc{Institut de Recherche Math\'ematique Avanc\'ee, Universit\'e de Strasbourg}\\
\noindent \textsc{7 Rue Ren\'e Descartes, 67084 Strasbourg CEDEX, France}\\
\noindent \emph{E-mail address:} \verb"robert.laterveer@math.unistra.fr"\\

\noindent \textsc{Institut de Recherche Math\'ematique Avanc\'ee, Universit\'e de Strasbourg}\\
\noindent \textsc{7 Rue Ren\'e Descartes, 67084 Strasbourg CEDEX, France}\\
\noindent \emph{E-mail address:} \verb"myriam.ounaies@unistra.fr"

\end{document}